\documentclass[a4paper,12pt]{amsart}
\usepackage[utf8]{inputenc}
\usepackage{amsmath,amssymb,amsfonts,amsthm,graphicx}

\newtheorem{theorem}{Theorem}[section]

\newtheorem{corollary}[theorem]{Corollary}
\numberwithin{equation}{section}

\theoremstyle{remark}

\newcommand{\Ric}{\mathop{\mathrm{Ric}}}

\title{On the Heat Kernel under the Ricci Flow Coupled with the Harmonic Map Flow}

\author{Mihai B\u{a}ile\c{s}teanu}
\thanks{Department of Mathematics, University of Rochester, 801 Hylan Bld, Rochester, NY 14627, USA \texttt{mbailest@z.rochester.edu} \\}

\begin{document}

\begin{abstract}
We estimate the heat kernel on a closed Riemannian manifold $M$, with $dim(M)\geq 3$, evolving under the Ricci-harmonic map flow and the result depends on some constants arising from a Sobolev imbedding theorem. In a special case, when the scalar curvature satisfies a certain natural inequality, we obtain, as a corollary, a bound similar to the one known for the fixed metric case. 
\end{abstract}

\maketitle

\section{Introduction}

We consider a coupled system of the Ricci flow, on a closed manifold $M$, with the harmonic map flow of a map $\phi$ from $M$ to some closed target manifold $N$: \begin{equation}\label{RH_flow_intro}
\begin{cases}\frac\partial{\partial t} g(x,t)=-2\Ric(x,t)+2\alpha(t)\nabla\phi(x,t)\otimes\nabla\phi(x,t)\\
\frac{\partial}{\partial t}\phi(x,t)=\tau_g\phi(x,t)
\end{cases}
 \end{equation} 
where $\alpha$ is a positive coupling time-dependent function. $\tau_g\phi$ denotes the tension field of the map $\phi$ with respect to the metric $g(t)$. We will refer to it as the $(RH)_\alpha$ flow, for short.  It is interesting that it may be less singular than both the Ricci flow (to which it reduces when $\alpha(t)=0$) and the harmonic map flow.  

The flow was introduced in \cite{RM12}, where the author studied several of its aspects: short time and long time existence, energy and entropy functionals, existence of singularities, a local non-collapsing property etc. The motivation for studying it was the work of List (\cite{BL08}), where the case of $\phi$ being a scalar function and $\alpha=2$ was analyzed. List found this form of the flow while searching for a connection between the classical Ricci flow and the static Einstein vacuum equations and discovered that this flow was equivalent to the gradient flow of an entropy functional, whose stationary points are solutions to the static Einstein vacuum equations. 

Another context in which the $(RH)\alpha$ flow arises is that of Ricci flow on warped product spaces. Given a warped product metric $g_M=g_N+e^{2\phi}g_F$ on a manifold $M=N\times F$ ($\phi\in C^{\infty}(N)$ ), the Ricci flow equation $\frac{\partial g_M}{\partial t}=-2\Ric_M$ leads to the following equations on each component: \[\begin{cases}
\frac{\partial g_N}{\partial t}=-2\Ric_N+2m\ d\phi\otimes\ d\phi \\
\frac{\partial \phi}{\partial t}=\triangle \phi-\mu e^{-2\phi}                                                                                                                                                                                                                                                                                                                                                                                                                                                                                                                                                                                                                                                       \end{cases}\] 
if the fibers $F$ are $m$-dimensional and $\mu$-Einstein. This is a modified version of the $(RH)_\alpha$, where the target manifold is one dimensional, and has been studied by M. B. Williams in \cite{MBW13} and by H. Tran in \cite{HT12} (when $\mu=0$).  

A natural direction in the study of the $(RH)_\alpha$ flow (which follows from the well-known analysis of the Ricci flow) is focusing on the behavior of the heat kernel. Generally, the study of the heat kernel has proven to be an interesting problem in its own right, but it has also given birth to many useful applications, for example in the study of the Ricci flow (the case with surgeries) or in the study of partial differential equations. The reason for that is that one can obtain information about the underlying structure from the behavior of the heat kernel.

For example, bounds on the heat kernel were found in relation to Harnack inequalities. In their seminal paper \cite{PLSTY86}, P. Li and S.-T. Yau obtained gradient estimates for positive solutions to the heat equation on closed manifolds with bounded Ricci curvature. Further they obtained Harnack inequalities and upper and lower bounds on the heat kernel. Their results were generalized to manifolds with non-convex boundary  by J. Wang in \cite{JW97}.

In recent years, there has been a growing interest in the study of the heat kernel on manifolds with time-dependent metrics. The validity of this question was settled by C. Guenther, in \cite{CG02}, where she studied the fundamental solution of a more general linear parabolic operator $L(u)=(\Delta-\frac{\partial}{\partial t}-h)u$, on compact n-dimensional manifolds with time dependent metric ($h$ being a smooth space-time function). All the expected properties of the heat kernel were shown to hold in this situation too: uniqueness, positivity, the adjoint property and the semigroup property. 

Moreover, a very important application of bounds on the heat kernel was shown by G. Perelman in \cite{GP02}, where he used a differential Li-Yau-Hamilton type inequalities for the fundamental solution of the conjugate heat equation $\Delta u+ \frac{\partial}{\partial t}-Ru=0$ under the Ricci flow to prove the pseudolocality theorem.  

Recently, X. Cao and Q. Zhang studied in \cite{XCQZ10} the Type I singularity model of the Ricci and proved that in the limit one obtains a gradient shrinking Ricci soliton by means of upper and lower bounds of the fundamental solution of the conjugate heat equation. 

The author obtained in \cite{MB12} a bound on the heat kernel under the Ricci flow, depending on the best constants in a Sobolev imbedding theorem. The techniques used were inspired by the work of Q. Zhang in \cite {QZ06}, where the conjugate heat equation introduced by Perelman (after a time reversal) is studied. 

In the present paper, we obtain a bound on the heat kernel under the $(RH)_\alpha$ flow and show, as a collorary, that when the manifold starts the flow with the scalar curvature satisfying $R>\alpha(0)|\nabla\phi(0)|$ at any point, then the bound becomes similar to the fixed metric case. This result, in fact, improves the result in \cite{MB12}, since in this case the constants are sharper. The result is based, as before, on a Sobolev imbedding theorem and there are no curvature assumptions. 

The main result is stated as follows:  
\begin{theorem}\label{theorem}
Let $M^n$ and $N^m$ be two closed Riemannian manifolds, with $n\geq 3$ and let $(g(t),\phi(t)), t\in[0,T]$ be a solution to the $(RH)_\alpha$ flow~\eqref{RH_flow_intro}, with $\alpha(t)$ a non-increasing positive function. Let $G(x,t;y,s)$ be the heat kernel, i.e.  fundamental solution for the heat equation $u_t=\bigtriangleup u$. Then there exists a positive number $C_n$, which depends only on the dimension $n$ of the manifold such that:

\begin{align*}
G(x,t;y,s)\leq &
\frac{C_n}{\left( \int\limits_{s}^{\frac{s+t}{2}}\left(\frac{m_0-c_n\tau}{m_0}\right)^{-2} \frac{e^{\frac{2}{n}H(\tau)}}{A(\tau)} \ d\tau\right)^{\frac{n}{4}} \left(\int\limits_{\frac{s+t}{2}}^{t}
\frac{e^{-\frac{2}{n}H(\tau)}}{A(\tau)} \ d\tau\right)^{\frac{n}{4}}} 
\end{align*}
for $0\leq s<t\leq T$; here $H(t)=\int\limits_{s}^{t}\left[\frac{B(\tau)}{A(\tau)}-\frac{3}{4}\cdot\frac{1}{m_0-c_n\tau}\right] d\tau $,  $1/m_0=\inf_{t=0}S$ - the infimum of $S=R-\alpha|\nabla\phi|^2$, taken at time $0$, and $A(t)$ and $B(t)$ are positive functions, which depend on the best constant in the Sobolev imbedding theorem stated above.
\end{theorem}

Notice that there are no curvature assumptions, but $B(t)$ will depend on the lower bound of the Ricci curvature and the derivatives of the curvature tensor at the initial time, as it will follow implicitly from Theorem \ref{thm_Aubin} presented below.  

The estimate may not seem natural, but in a special case, when the scalar curvature $R(x,0)>\alpha(0)|\nabla\phi(x,0)|^2$, one obtains a bound similar to the the fixed metric case. Recall that J. Wang obtained in \cite{JW97} that the heat kernel on an $n$-dimensional compact Riemannian manifold $M$, with fixed metric, is bounded from above by $N(S)(t-s)^{-n/2}$, where $N(S)$ is the Neumann Sobolev constant of $M$, coming from a Sobolev imbedding theorem. Our corollary exibits a similar bound:   

\begin{corollary}
Under the same assumptions as in theorem (\ref{theorem}), together with the condition that $R(x,0)>\alpha(0)|\nabla\phi(x,0)|^2$, there exists a positive number $\tilde{C}_n$, which depends only on the dimension $n$ of the manifold and on the best constant in the Sobolev imbedding theorem in $\mathbb{R}^n$, such that:
\begin{align*}
G(x,t;y,s)\leq \tilde{C}_n\cdot\frac{1}{(t-s)^{\frac{n}{2}}} \hspace{1cm} \text{ for } 0\leq s<t\leq T
\end{align*}
\end{corollary}

The exact expression of $\tilde{C}_n$ is $\left(\frac{4K(n,2)}{n}\right)^{\frac{n}{2}}$, where $K(n,2)$ is the best constant in the Sobolev imbedding in $\mathbb{R}^n$.

The paper is structured as follows: in section \ref{doi} we will introduce the notation, state the Sobolev imbedding theorems used, while the proofs of the theorem and of the corollary are presented in section \ref{patru}. 

The author would like to acknowledge prof. Xiaodong Cao for continuous support and valuable advice.  

\section{Background and notation}\label{doi}

Let $(M^n,g)$ and $(N^m, \gamma)$ two $n$-dimensional and $m$-dimensional ($n\geq 3$), respectively, manifolds without boundary, which are compact, connected, oriented and smooth. We also let $g(t)$ be a family of Riemannian metrics on $M$, while $\phi(t)$ a family of smooth maps between $M$ and $N$.    

For $T>0$, let $(g(t),\phi(t)), t\in[0,T]$ be a solution to the coupled system of Ricci flow and harmonic map flow, i.e. the $(RH)_\alpha$ flow,  with coupling time-dependent constant $\alpha(t)$: \begin{equation}\label{RH_flow}
\begin{cases}\frac\partial{\partial t} g(x,t)=-2\Ric(x,t)+2\alpha(t)\nabla\phi(x,t)\otimes\nabla\phi(x,t)\\
\frac{\partial}{\partial t}\phi(x,t)=\tau_g\phi(x,t)
\end{cases}
 \end{equation} 
 
Note that $\nabla\phi(x,t)\otimes\nabla\phi(x,t)$ has the following expression in local coordinates: $(\nabla\phi\otimes\nabla\phi)_{ij}=\nabla_i\phi^\lambda\nabla_j\phi^\lambda$.

We restrict our attention to a time interval $[0,T]$ where a solution to this system exists (R. M\"uller showed in \cite{RM12} that the flow has short time existence, and we will not worry about the behavior at the blow-up time). We are interested in obtaining bounds on the heat kernel $G(x,t;y,s)$, which is the fundamental solution of the heat equation
\begin{align}\label{heateqn}
\left(\Delta-\frac\partial{\partial t}\right)u(x,t)=0,\qquad x\in M,~t\in[0,T].
\end{align}
Such heat kernel does indeed exist and it is well defined, as it was shown in \cite{CG02} by C. Guenther, who studied the fundamental solution of the linear parabolic operator $L(u)=(\Delta-\frac{\partial}{\partial t}-h)u$, on compact n-dimensional manifolds with time dependent metric, where $h$ is a smooth space-time function. She proved the uniqueness, positivity, the adjoint property and the semigroup property of this operator, which thus behaves like the usual heat kernel. As a particular case ($h=0$), she obtained the existence and properties of the heat kernel under any flow of the metric.

We denote the heat kernel by $G(x,t;y,s)$ and it is a smooth function $G(x,t;y,s):M\times[0,T]\times M\times[0,T]\to\mathbb{R}$, with $s<t$, which satisfies two properties: 
\begin{enumerate}
\item[(i)] $L(G)=0$ in $(x,t)$ for $(x,t)\neq (y,s)$ 
\item[(ii)] $\lim_{t\to s}G(.,t;y,s)=\delta_y$ for every $y$, where $\delta_y$ is the Dirac delta function.
\end{enumerate} 

By definition, $G$ satisfies the heat equation in the $(x,t)$ coordinates \[\Delta_x G(x,t;y,s)-\partial_t G(x,t;y,s)=0\] whereas in the $(y,s)$ it satisfies the adjoint (conjugate) heat equation \[\Delta_y G(x,t;y,s)+\partial_s G(x,t;y,s)- [R(y,s)-\alpha|\nabla\phi|^2] G(x,t;y,s)=0\] (see \cite{RM12} for a proof of this fact), where $R(y,s)$ is the scalar curvature, measured with respect to the metric $g(s)$.

Following the notation in \cite {RM12}, it will be easier to introduce these quantities: 
\begin{align*}
 \mathcal{S}&:= \Ric-\alpha\nabla\phi\otimes\nabla\phi\\
 S_{ij}&:=R_{ij}-\alpha\nabla_i\phi\nabla_j\phi\\
 S&:=R-\alpha|\nabla\phi|^2
\end{align*}

Next, we will present the Sobolev inequalities that form the basis of our exploration. 

Sharpening a result by T. Aubin (\cite{TA76}), E. Hebey proved in (\cite{EHMV96}) the following : 
\begin{theorem}\label{thm_Aubin}
Let $M^n$ be a smooth compact Riemannian manifold of dimension $n$. Then there exists a constant $B$ such that for any $\psi\in W^{1,2}(M)$ (the Sobolev space of weakly differentiable functions) :
\[||\psi||_p^2\leq K(n,2)^2||\nabla\psi||_2^2+B||\psi||_{2}^2 \hspace{0.5cm}.\]

Here $K(n,2)$ is the best constant in the Sobolev imbedding (inequality) in $\mathbb{R}^n$ and $p=(2n)/(n-2)$. $B$ depends on the lower bound of the Ricci curvature and the derivatives of the curvature tensor. 
\end{theorem}

Note that Hebey's result was shown for complete manifolds, and in that situation $B$ depends on the injectivity radius. However, we are interested in compact manifolds, so $B$ will not depend on the injectivity radius. 

Later, along the Ricci flow, Q. Zhang proved the following uniform Sobolev inequality in \cite{QZ09}:

\begin{theorem}\label{thm_Zhang}
Let $M^n$ be a compact Riemannian manifold, with $n\geq 3$ and let $\big(M,g(t)\big)_{t\in[0,T]}$ be a solution to the Ricci flow $\frac{\partial g}{\partial t}=-2\Ric$. Let $A$ and $B$ be positive numbers such that for $(M,g(0))$ the following Sobolev inequality holds: for any $v\in W^{1,2}(M,g(0))$,
\[\left(\int_M|v|^{\frac{2n}{n-2}}\ d\mu(g(0))\right)^{\frac{n-2}{n}}\leq A\int_M|\nabla v|^2\ d\mu(g(0))+B\int_Mv^2\ d\mu(g(0))\]

Then there exist positive functions $A(t)$, $B(t)$ depending only on the initial metric $g(0)$ in terms of $A$ and $B$, and $t$ such that, for all $v\in W^{1,2}(M,g(t))$, $t>0$, the following holds
 
\[\left(\int_M|v|^{\frac{2n}{n-2}}\ d\mu(g(t))\right)^{\frac{n-2}{n}}\leq A(t)\int_M\left(|\nabla v|^2+\frac{1}{4}Rv^2\right)\ d\mu(g(t))+B(t)\int_Mv^2\ d\mu(g(t))\]

Here $R$ is the scalar curvature with respect to $g(t)$. Moreover, if $R(x,0)>0$, then $A(t)$ is independent of $t$ and $B(t)=0$.
\end{theorem}

The proof of this theorem relies on the analysis of $\lambda_0$, which is the first eigenvalue of Perelman's $\mathcal{F}$-entropy, i.e.:
\[\lambda_0=\inf\limits_{||v||_2=1}\int_M(4|\nabla v|^2+Rv^2)\ d\mu(g(0)).\]

Recently it has been proven that the same theorem holds for the Ricci flow coupled with the harmonic map flow (see \cite{LXGWK13}), where the analysis is now based on 
\[\lambda_0^\alpha=\inf\limits_{||v||_2=1}\int_M(4|\nabla v|^2+Sv^2)\ d\mu(g(0)),\]
where $S=R-\alpha|\nabla\phi|^2$. In the new setting, at time $t$ there are positive functions $A(t)$ and $B(t)$ such that   
\[\left(\int_M|v|^{\frac{2n}{n-2}}\ d\mu(g(t))\right)^{\frac{n-2}{n}}\leq A(t)\int_M\left(|\nabla v|^2+\frac{1}{4}Sv^2\right)\ d\mu(g(t))+B(t)\int_Mv^2\ d\mu(g(t))\]

Moreover, if $\lambda_0^\alpha>0$, which is automatically satisfied if at the initial time $R(0)>\alpha (0)|\nabla\phi(0)|$, then $A(t)$ is a constant and $B(t)=0$.  

Note that R. M\"uller introduced in \cite{RM12} the $\mathcal{F}_{\alpha}$ and the $\mathcal{W}_\alpha$ functionals, which are the natural analogues of the $\mathcal{F}$ and the $\mathcal{W}$ functionals for the Ricci flow, introduced by Perelman.

\section{Proof}\label{patru}

%\begin{proof} 

We start the proof by assuming, without loss of generality, that $s=0$. By the semigroup property of the heat kernel \cite[Theorem 2.6]{CG02} and the Cauchy-Bunyakovsky-Schwarz inequality we have that:
\begin{align*}
G(x,t;y,0)& =\int_M G\left(x,t;z,\frac{t}{2}\right)G\left(z,\frac{t}{2};y,0\right)d\mu\left(z,\frac{t}{2}\right) \\
          & \leq \left[\int_M G^2\left(x,t;z,\frac{t}{2}\right)d\mu\left(z,\frac{t}{2}\right) \right]^{1/2}\left[\int_M G^2\left(z,\frac{t}{2};y,0\right)d\mu\left(z,\frac{t}{2}
\right)\right]^{1/2}
\end{align*}

The key of the proof consists in determining upper bounds for the following two quantities:
\begin{align*}
\alpha(t)=\int_M G^2(x,t;y,s)d\mu(x,t)& \text{ (for $s$ fixed)}\\
\beta(s)=\int_M G^2(x,t;y,s)d\mu(y,s)& \text{ (for $t$ fixed)}
\end{align*}

We will find an ordinary differential inequality for each of the two.

We first deduce a bound on $\alpha(t)$, by finding an inequality involving $\alpha'(t)$ and $\alpha(t)$. Note that we will treat $G$ as being a function of $(x,t)$, the $(y,s)$ part is fixed. 

Since $\frac{d}{dt}(d\mu)=-Sd\mu=(-R+\alpha|\nabla\phi|^2)d\mu$ (due to the Ricci-harmonic flow), one has:

\begin{align}\label{pprim}
\alpha'(t) & =2\int_M G\cdot G_t\ d\mu(x,t)-\int_M G^2 (R-\alpha|\nabla\phi|^2)\ d\mu(x,t) \notag\\
 & = 2\int_M G\cdot(\Delta G)\ d\mu(x,t)-\int_M G^2 (R-\alpha|\nabla\phi|^2)\ d\mu(x,t)\notag \\
      &= -2\int_M|\nabla G|^2\ d\mu-\int_M G^2(R-\alpha|\nabla\phi|^2)\ d\mu \notag \\ 
      & \leq -\int_M[|\nabla G|^2+(R-\alpha|\nabla\phi|^2) G^2]\ d\mu(x,t)
\end{align}

Estimating $\int_M|\nabla G|^2d\mu$ will make use of the Sobolev imbedding theorem, which gives a relation between $\int_M|\nabla G|^2d\mu$ and $\int_M G^2d\mu$, and the H\"older inequality to bound the term involving $G^{2n/(n-2)}$:
\begin{align}\label{Holder}
\int_M G^2\ d\mu(x,t)\leq \left[\int_M G^{\frac{2n}{n-2}}\ d\mu(x,t) \right]^{\frac{n-2}{n+2}} \left[\int_M G\ d\mu(x,t) \right]^{{4\over n+2}}
\end{align}

By theorem (\ref{thm_Aubin}) one gets that at time $t=0$, the following inequality holds for any $v\in W^{1,2}(M,g(0))$ (hence also for $G(x,t:y,s)$, which is smooth) and for some $B>0$:

\[\left(\int_M|v|^{\frac{2n}{n-2}}\ d\mu(g(0))\right)^{\frac{n-2}{n}}\leq K(n,2)^2\int_M|\nabla v|^2\ d\mu(g(0))+B\int_Mv^2\ d\mu(g(0))\]

Then by Theorem \ref{thm_Zhang} (but applied for the Ricci-harmonic flow) it follows that at any time $t\in (0,T]$ and for all $v\in W^{1,2}(M,g(t))$:

\[\left(\int_M|v|^{\frac{2n}{n-2}}\ d\mu(g(t))\right)^{\frac{n-2}{n}}\leq A(t)\int_M\left(|\nabla v|^2+\frac{1}{4}(R-\alpha|\nabla\phi|^2)v^2\right)\ d\mu(g(t))+B(t)\int_Mv^2\ d\mu(g(t))\]

where $A(t)$ is a positive function depending on $g(0)$ and $K(n,2)^2$, while  $B(t)$ is also a positive function, depending on $B$, which in turn depends on the initial Ricci curvature on $M$ and on the derivatives of the curvatures on $M$ at time $0$.

Applying the above for the heat kernel, one can relate the RHS of (\ref{Holder}) to the Sobolev inequality:
\begin{gather}\label{Sobolev1}
\int_M G^2\ d\mu(x,t)  \leq \left[\int_M G^{\frac{2n}{n-2}}\ d\mu(x,t) \right]^{\frac{n-2}{n+2}} \left[\int_M G\ d\mu(x,t) \right]^{{4\over n+2}} \notag\\
\leq \left[ A(t)\int_M\left(|\nabla G|^2+\frac{1}{4}(R-\alpha|\nabla\phi|^2)G^2\right)\ d\mu(x,t)+B(t)\int_MG^2\ d\mu(x,t)\right]^{\frac{n}{n+2}}\left[\int_M G\ d\mu(x,t)\right]^{\frac{4}{n+2}}
\end{gather}

We will now focus our attention to $J(t):=\int_M G(x,t;y,s)\ d\mu(x,t)$. By the definition of the fundamental solution $\int_M G(x,t;y,s)d\mu(y,s)=1$, but that's not true if one integrates in $(x,t)$. Our goal will be to obtain a differential inequality for $J(t)$, from which a bound will be found. 

\begin{align*}
J'(t)&=\int_M G_t(x,t;y,s)\ d\mu(x,t)+\int_M G(x,t;y,s)\frac{d}{dt}\ d\mu(x,t)=\int_M\Delta_xG(x,t;y,s)\ d\mu(x,t)\\
     &-\int_M G(x,t;y,s)S(x,t)\ d\mu(x,t)=-\int_M G(x,t;y,s)S(x,t)\ d\mu(x,t)
\end{align*}
the first term being $0$, as $M$ is a compact manifold, without boundary.

According to Theorem 4.4 from \cite{RM12}, $S$ satisfies the following equation:
\begin{align*}
\frac{\partial S}{\partial t}=\bigtriangleup S+2|S_{ij}|^2+2\alpha|\tau_g\phi|^2-\frac{\partial\alpha}{\partial t}|\nabla\phi|^2
\end{align*}
But $|S_{ij}|^2\geq \frac{1}{n}S^2$ (this is true for any 2-tensor) and $\alpha(t)$ is a positive function, non-increasing in time, so one obtains: 
\[\frac{\partial S}{\partial t}-\bigtriangleup S-\frac{2}{n}S^2\geq 0\]

Since the solutions of the ODE $\frac{d\rho}{dt}=\frac{2}{n}\rho^2$
are $\rho(t)=\frac{n}{n\rho(0)^{-1}-2t}$, by the maximum principle
we get a bound on $S$, for $s\leq \tau\leq t$:
\begin{align*}
S(z,\tau)\geq \frac{n}{n(\inf_{t=0}S)^{-1}-2\tau}=\frac{1}{(\inf_{t=0}S)^{-1}-\frac{2}{n}\tau}:=\frac{1}{m_0-c_n\tau}
\end{align*}
(here and later, if $\inf_{t=0} S\geq 0$, then the above is regarded
as zero).

\bigskip
Using this lower bound for $S$ (for $\tau\in(s,t]$), we get:
\begin{align*}
J'(\tau)\leq -\frac{1}{m_0-c_n\tau}J(\tau)
\end{align*}
After integrating the above from $s$ to $t$, while noting that by $J(s)$ one understands: 
\[J(s)=\lim\limits_{t\to s}\int_M G(x,t;y,s)\ d\mu(x,t)=\int_M \lim\limits_{t\to s} G(x,t;y,s)\ d\mu(x,t)=\int_M \delta_{y}(x)\ d\mu(x,s)=1\]
one obtains:
\begin{align*}
J(t)\leq \left(\frac{m_0-c_nt}{m_0-c_ns}\right)^{\frac{n}{2}}:=(\chi_{t,s})^{\frac{n}{2}}
\end{align*}

Hence $\int_M G(x,t;y,s)\ d\mu(x,t)\leq (\chi_{t,s}))^{\frac{n}{2}}$ and (\ref{Sobolev1}) becomes:

\begin{align*}
\int_M G^2d\mu(x,t)\leq \left[ A(t)\int_M\left(|\nabla G|^2+\frac{1}{4}SG^2\right)\ d\mu(x,t)+B(t)\int_MG^2\ d\mu(x,t)\right]^{\frac{n}{n+2}}\left(\chi_{t,s}\right)^{\frac{2n}{n+2}}
\end{align*}

From this it follows immediately that:

\begin{align*}
\int_M |\nabla G|^2\ d\mu(x,t)\geq \frac{1}{\chi^2_{t,s}A(t)}\left[\int_M G^2\ d\mu(x,t)\right]^{\frac{n+2}{n}} -\frac{B(t)}{A(t)}\int_M G^2\ d\mu(x,t)-\frac{1}{4}\int SG^2\ d\mu(x,t)
\end{align*}

Combining this with the inequality from (\ref{pprim}), one obtains the following differential inequality for $\alpha(t)$:
\begin{align*}
\alpha'(t)\leq -\frac{1}{\chi^2_{t,s}A(t)} \alpha(t)^{{n+2\over n}}+\frac{B(t)}{A(t)}\alpha(t)-\frac{3}{4}\int SG^2d\mu(x,t)
\end{align*}

Note that the above is true for any $\tau\in(s,t]$. For the following computation, we will consider $t$ fixed as well.  Recall that for $\tau\in (s,t]$, $S(\cdot, \tau)\geq\frac{1}{m_0-c_n\tau}$. Denoting with:
\begin{align*}
h(\tau):=\frac{B(\tau)}{A(\tau)}-\frac{3}{4}\cdot\frac{1}{m_0-c_n\tau}
\end{align*}
we get:
\begin{align*}
\alpha'(\tau)\leq -\frac{1}{\chi^{2}_{\tau,s}A(\tau)} \alpha(\tau)^{{n+2\over n}}+h(\tau)\alpha(\tau)
\end{align*}

Let $H(\tau)$ be an antiderivative of $h(\tau)$. By the integrating factor method, one finds:
\begin{align*}
(e^{-H(\tau)}\alpha(\tau))'\leq  -\frac{1}{\chi^{2}(\tau)A(\tau)}(e^{-H(\tau)}\alpha(\tau))^{\frac{n+2}{n}}e^{\frac{2}{n}H(\tau)}
\end{align*}
Since the above is true for any $\tau\in(s,t]$, by integrating from $s$ to $t$ and taking into account that 
\[\lim\limits_{\tau\searrow s}\alpha(\tau)=\int_M\lim\limits_{\tau\searrow s}G^2(x,\tau;y,s)\ d\mu(x,\tau)=\int_M\delta^2_y(x)\ d\mu(x,s)=0\]
one obtains the first necessary bound:
\begin{align*}
\alpha(t)\leq \frac{C_ne^{H(t)}}{\left(\int\limits_{s}^{t}\frac{e^{\frac{2}{n}H(\tau)}}{\chi^{2}(\tau)A(\tau)} d\tau\right)^{\frac{n}{2}}}
\end{align*}
where $C_n=\left(\frac{2}{n}\right)^{\frac{n}{2}}$.

\bigskip
The next step is to estimate $\beta(s)=\int_M G^2(x,t;y,s)\ d\mu(y,s)$, for which the computation is different, due to the assymetry of the
equation. As stated above, the second entries of $G$ satisfy the conjugated equation:
\begin{align*}
\Delta_y G(x,t;y,s)+\partial_s G(x,t;y,s)- S G(x,t;y,s)=0
\end{align*}

Proceeding just as in the $\alpha(t)$ case, we get the following:
\begin{align*}
\beta'(s)& =2\int_M GG_s\ d\mu(y,s)-\int_M SG^2\ d\mu(y,s)\\
         &=2\int_M G(-\Delta G +SG)\ d\mu(y,s)-\int_M SG^2\ d\mu(y,s)\\
       & =-2\int_M G(\Delta G)\ d\mu(y,s)+\int_M SG^2\ d\mu(y,s)\\
       & =2\int_M |\nabla G|^2\ d\mu(y,s)+\int_M SG^2\ d\mu(y,s) \\
     & \geq \int_M |\nabla G|^2\ d\mu(y,s)+\int_M SG^2\ d\mu(y,s)
\end{align*}

Hence
\begin{align*}
\beta'(s)\geq \int_M (|\nabla G|^2+SG^2)\ d\mu(y,s)
\end{align*}

But this time, by the property of the heat kernel:
\begin{align*}
\tilde{J}(s):=\int_M G(x,t;y,s)\ d\mu(y,s)=1
\end{align*}

so by applying H\"older (as for $\alpha(t)$) and relating it to the Sobolev inequality, we get:
\begin{align*}
\int_M G^2\ d\mu(y,s)&\leq \left[ A(s)\int_M\left(|\nabla G|^2+\frac{1}{4}SG^2\right)\ d\mu(y,s)+B(s)\int_MG^2\ d\mu(y,s)\right]^{\frac{n}{n+2}}\left[\int_M G\ d\mu(y,s)\right]^{\frac{4}{n+2}} \\
 & =\left[ A(s)\int_M\left(|\nabla G|^2+\frac{1}{4}SG^2\right)\ d\mu(y,s)+B(s)\int_MG^2\ d\mu(y,s)\right]^{\frac{n}{n+2}}
\end{align*}

Following the same steps as for $\alpha(t)$, one finds
\begin{align*}
\beta'(s)\geq \frac{1}{A(s)} \beta(s)^{{n+2\over n}}-h(s)\beta(s)
\end{align*}
($h(s)$ denotes, as before, $\frac{B(s)}{A(s)}-\frac{3}{4}\cdot\frac{1}{m_0-c_ns}$)

The above is true for any $\tau\in[s,t)$. We will apply again the integrating factor method, with $H(\tau)$ being the same antiderivative of $h(\tau)$ as above. 
For $\tau\in [s,t)$, the following holds:
\begin{align*}
(e^{H(\tau)}\beta(\tau))'\geq \frac{1}{A(\tau)}(e^{H(\tau)}\beta(\tau))^{{n+2\over n}}e^{-\frac{2}{n}H(\tau)}
\end{align*}

Integrating between $s$ and $t$, and taking into account that
\[ \lim\limits_{\tau\nearrow t}\beta(\tau)=\int_M\lim\limits_{\tau\nearrow t}G^2(x,t;y,\tau)\ d\mu(y,\tau)=\int_M\delta^2_y(x)\ d\mu(y,t)=0\] 
we get the second desired bound:
\begin{align*}
\beta(s)\leq \frac{C_ne^{-H(s)}}{\left(\int\limits_{s}^{t} \frac{e^{-\frac{2}{n}H(\tau)}}{A(\tau)}\ d\tau\right)^{n/2}}
\end{align*}

From the estimates of $\alpha$ and $\beta$ we get the following:

\begin{align*}
\alpha\left(\frac{t}{2}\right)=\int_M G^2\left(z,\frac{t}{2};y,0\right)\ d\mu\left(z,t/2\right)\leq \frac{C_ne^{H(t/2)}}{\left( \int\limits_{0}^{t/2}\left(\frac{m_0-c_n\tau}{m_0}\right)^{-2} \frac{e^{\frac{2}{n}H(\tau)}}{A(\tau)} \ d\tau\right)^{\frac{n}{2}}}
\end{align*}
\begin{align*}
\beta\left(\frac{t}{2}\right)=\int_M G^2\left(x,t;z,\frac{t}{2}\right)\ d\mu\left(z,\frac{t}{2}\right)\leq \frac{C_ne^{-H(t/2)}}{\left(\int\limits_{t/2}^{t}
\frac{e^{-\frac{2}{n}H(\tau)}}{A(\tau)} \ d\tau\right)^{n/2}}
\end{align*}

Here, we may choose $H(t/2)=\int\limits_{0}^{t/2}\left[\frac{B(\tau)}{A(\tau)}-\frac{3}{4}\cdot\frac{1}{m_0-c_n\tau}\right]\ d\tau$, since the relation is true for any antiderivative of $h(\tau)=\frac{B(\tau)}{A(\tau)}-\frac{3}{4}\cdot\frac{1}{m_0-c_n\tau}$.

The conclusion follows from multiplying the relations above.

%\end{proof}

\subsection{Proof of the collorary}

In the special case when $S(x,0)>0$, then $S(x,t)>0$ for all $t>0$, so it follows that $J'(\tau)\leq 0$. $J(\tau)$ is thus decreasing, so $J(\tau)\leq J(s)=1$, which leads to the differential inequality for $\alpha(t)$ to be: 

\begin{align*}
\alpha'(t)\leq -\frac{1}{A(t)} \alpha(t)^{{n+2\over n}}+\frac{B(t)}{A(t)}\alpha(t)
\end{align*}

And from this the bound for $\alpha(t)$ becomes:

\begin{align*}
\alpha(t)\leq \frac{C_ne^{H(t)}}{\left(\int\limits_{s}^{t}\frac{e^{\frac{2}{n}H(\tau)}}{A(\tau)} d\tau\right)^{\frac{n}{2}}}
\end{align*}

where $H(\tau)$ is the antiderivative of $\frac{B(\tau)}{A(\tau)}$ such that $H(s)\neq 0$ and $H(t)\neq 0$.

Similarly, one obtaines for $\beta(s)$:

\begin{align*}
\beta'(s)\geq \frac{1}{A(s)} \beta(s)^{{n+2\over n}}-\frac{B(s)}{A(s)}\alpha(s)
\end{align*}

and from this: 

\begin{align*}
\beta(s)\leq \frac{C_ne^{-H(s)}}{\left(\int\limits_{s}^{t} \frac{e^{-\frac{2}{n}H(\tau)}}{A(\tau)}\ d\tau\right)^{n/2}}
\end{align*}

where $H(\tau)$ is the same antiderivative of $\frac{B(\tau)}{A(\tau)}$ as above. 

By (\ref{thm_Zhang}), in the case of $S(x,0)>0$, the function $A(t)$ is a constant, while $B(t)=0$. Recall that $A(t)=A(0)$ is in fact $K(n,2)$, where $K(n,2)$ is the best constant in the Sobolev imbedding.

One has then that $H(t)=\frac{B}{A}t=0$. Using this, we get:
\begin{align*}
G(x,t;y,s)\leq &
\frac{C_n}{\left(\int\limits_{s}^{\frac{s+t}{2}}\frac{1}{A(0)} \ d\tau\right)^{\frac{n}{4}} \left(\int\limits_{\frac{s+t}{2}}^{t}
\frac{1}{A(0)} \ d\tau\right)^{\frac{n}{4}}}=\frac{C_n}{\left[\left(\frac{t-s}{2A}\right)^2\right]^{\frac{n}{4}}} = \frac{\tilde{C}_n}{(t-s)^{\frac{n}{2}}}
\end{align*}

where $\tilde{C}_n=C_n\cdot (2A)^{\frac{n}{2}}=\left(\frac{4K(n,2)}{n}\right)^{\frac{n}{2}}$.

This proves the desired corollary.

\subsection{Final remarks}

We have shown that under the Ricci-harmonic $(RH)_\alpha$ flow, if one starts with the condition that $R>\alpha|\nabla\phi|^2$ (which is preserved under the flow) the heat kernel behaves similarlty to the case when the metric is static. It is worthwile to investigate if one can prove Gaussian bounds under a more general setting.

\bibliographystyle{unsrt}
\bibliography{Heat-kernel-estimate}

\end{document}